\documentclass{IEEEtran} 
\usepackage[latin1]{inputenc}
\usepackage{cite}
\usepackage{color}
\usepackage{amsmath,amsfonts,amssymb}

\DeclareMathOperator{\QFT}{QFT}
\usepackage{url,xspace}
\usepackage{graphicx}
\newtheorem{lemma}{Lemma}
\newtheorem{definition}{Definition}
\newtheorem{proposition}{Proposition}

\newtheorem{property}{Property}
\renewcommand{\v}[1]{\boldsymbol{#1}} 
\renewcommand{\i}{\ensuremath{\v{i}}\xspace}
\renewcommand{\j}{\ensuremath{\v{j}}\xspace}
\renewcommand{\k}{\ensuremath{\v{k}}\xspace}

\newcommand{\norm}[1]{\ensuremath{||#1||}\xspace}

\newcommand{\R}{\ensuremath{{\mathbb R}}}
\newcommand{\C}{\ensuremath{{\mathbb C}}}
\renewcommand{\H}{\ensuremath{{\mathbb H}}}
\newcommand*{\qconjugate}[1]{\overline{#1}}                
\newcommand*{\cconjugate}[1]{{#1^\star}}                   
\DeclareMathOperator{\sign}{sign}
\newcommand{\dt}{\mathrm{d}t}
\newcommand{\dnu}{\mathrm{d}\nu}
\newcommand{\dtau}{\mathrm{d}\tau}

\newcommand*{\hilbert}[1]{\ensuremath{\mathcal{H}\left[#1\right]}\xspace}
\newcommand*{\qhilbert}[2]{\ensuremath{\mathcal{H}_{#1}\left[#2\right]}\xspace}

\newcommand*{\fourier}{\ensuremath{\mathcal{F}}\xspace}
\newcommand*{\ffourier}[1]{\ensuremath{\mathcal{F}\left[#1\right]}\xspace} 
\newcommand*{\ifourier}[1]{\ensuremath{\mathcal{F}^{-1}\left[#1\right]}\xspace}
\newcommand*{\qfourier}[2]{\ensuremath{\mathcal{F}_{#1}\left[#2\right]}\xspace} 
\newcommand*{\iqfourier}[2]{\ensuremath{\mathcal{F}_{#1}^{-1}\left[#2\right]}\xspace}
\let\oldmu\mu\renewcommand{\mu}{\v{\oldmu}} 
\let\oldxi\xi\renewcommand{\xi}{\v{\oldxi}}

\begin{document}

\title{The hyperanalytic signal}

\author{Nicolas~Le~Bihan~and~Stephen~J.~Sangwine~\IEEEmembership{Senior~Member,~IEEE}%
\thanks{N. Le Bihan is with CNRS, GIPSA-Lab, Département Images et Signal,
        961 Rue de la Houille Blanche, Domaine Universitaire,
        BP 46, 38402 Saint Martin d'Hères cedex, France.}%
\thanks{S. J. Sangwine is with the School of Computer Science and Electronic Engineering,
        University of Essex, Colchester, CO4 3SQ, UK.}}

\maketitle

\begin{abstract}
The concept of the analytic signal is extended from the case
of a real signal with a complex analytic signal to a complex
signal with a hypercomplex analytic signal
(which we call a hyperanalytic signal)
The hyperanalytic signal may be interpreted as an ordered pair
of complex signals or as a quaternion signal.
The hyperanalytic signal contains a complex orthogonal signal
and we show how to obtain this by three methods: a pair of
classical Hilbert transforms; a complex Fourier transform;
and a quaternion Fourier transform.
It is shown how to derive from the hyperanalytic signal a
complex envelope and phase using a polar quaternion representation
previously introduced by the authors.
The complex modulation of a real sinusoidal carrier is shown
to generalize the modulation properties of the classical
analytic signal.
The paper extends the ideas of properness to deterministic
complex signals using the hyperanalytic signal.
A signal example is presented, with its orthogonal signal, and
its complex envelope and phase.
\end{abstract}

\section{Introduction}
\label{Introsection}

The analytic signal has been known since 1948 from
the work of Ville \cite{Ville:1948} and Gabor \cite{Gabor:1946}.
It can be simply described, even though its theoretical
ramifications are deep. Its use in non-stationary signal
analysis is routine and it has been used in numerous applications.
Simply put, given a real-valued signal $f(t)$, its analytic signal $a(t)$ is a complex
signal with real part equal to $f(t)$ and imaginary part orthogonal to $f(t)$.
The imaginary part is sometimes known as the quadrature signal -- in the case where
$f(t)$ is a sinusoid, the imaginary part of the analytic signal is in quadrature,
that is with a phase difference of $-\pi/2$.
The orthogonal signal is related to $f(t)$ by the Hilbert transform \cite{Hahn:1996,Hahn:Poularikas}.
The analytic signal has the interesting property that its modulus $|a(t)|$ is an
envelope of the signal $f(t)$.
The envelope is also known as the \emph{instantaneous amplitude}.
Thus if $f(t)$ is an amplitude-modulated sinusoid, the envelope $|a(t)|$, subject to
certain conditions, is the modulating signal. The argument of the analytic signal,
$\angle\,a(t)$ is known as the \emph{instantaneous phase}.
The analytic signal has a third interesting property: it has a one-sided Fourier transform.
Thus a simple algorithm for constructing the analytic signal (algebraically or numerically)
is to compute the Fourier transform of $f(t)$,
multiply the Fourier transform by a unit step which is zero for negative frequencies,
and then construct the inverse Fourier transform.

In this paper we extend the concept of the analytic signal from the case of a real signal
$f(t)$ with a complex analytic signal $a(t)$,
to a complex signal $z(t)$ with a \emph{hypercomplex} analytic signal $h(t)$,
which we call the \emph{hyperanalytic signal}.
Just as the classical complex analytic signal contains
both the original real signal (in the real part) and a real orthogonal signal
(in the imaginary part),
a hyperanalytic signal contains two complex signals:
the original signal and an orthogonal signal.
We have previously published partial results on this topic
\cite{SangwineLeBihan:2007,LeBihanSangwine:2008,LeBihanSangwine:2008a}.
In the present paper we develop a clear idea of how to generalise the classic
case of amplitude modulation to a complex signal, and show for the first time
that this leads to a correctly extended analytic signal concept in which the
(\emph{complex}) envelope and phase have clear interpretations.
We show how the orthogonal signal can be constructed using
three methods, yielding the same result:
\begin{enumerate}
\item a pair of classical Hilbert transforms operating independently on the
      real and imaginary parts of $z(t)$,
\item a complex Fourier transform pair operating on $z(t)$ as a whole, with
      positive and negative frequency coefficients rotated by $\pm\pi/2$
      using $i$ times a signum function,
\item a quaternion Fourier transform with suppression of negative frequencies
      in the Fourier domain.
\end{enumerate}
All three of these methods are based in one way or another on single-sided spectra,
but it is only in the third (quaternion) case that a single one-sided spectrum is involved.

The construction of an orthogonal signal alone would not constitute a full
generalisation of the classical analytic signal to the complex case: it is
also necessary to generalise the envelope and phase concepts, and this can
only be done by interpreting the original and orthogonal complex signals
as a pair.
In this paper we have only one way to do this: by representing the pair of
complex signals as a quaternion signal.
This arises naturally from the third method above for creating an orthogonal
signal, but also from the Cayley-Dickson construction of a quaternion as a
complex number with complex real and imaginary parts (with different
roots of $-1$ used in each of the two levels of complex number).

Extension of the analytic signal concept to 2D signals, that is images,
with real, complex or quaternion-valued pixels is of interest, but outside the scope of this paper.
Some work has been done on this, notably by Bülow, Felsberg, Sommer and Hahn \cite{Bulow:1999,BulowSommer:2001,FelsbergSommer:2001,10.1109/5.158601}.
The principal issue to be solved in the 2D case is to generalise the concept of
a single-sided spectrum. Hahn considered a single quadrant or orthant spectrum,
Sommer \textit{et al} considered a spectrum with support limited to half the
complex plane, not necessarily confined to two quadrants,
but still with real sample or pixel values.

Recently, Lilly and Olhede \cite{10.1109/TSP.2009.2031729} have published a paper
on bivariate analytic signal concepts without explicitly considering the complex
signal case which we cover here. Their approach is linked to a specific signal
model, the \emph{modulated elliptical signal}, which they illustrate with the
example of a drifting oceanographic float. The approach taken in the present paper
is more general and without reference to a specific signal model.

\section{Construction of an orthogonal signal using classical complex transforms}
\label{Hypersection}

\label{sec:hyperanalytic}
We first review and examine the notion of orthogonality in the
case of complex signals before we consider how to construct a signal
orthogonal to a given complex signal.
\begin{definition}
\label{def:orthogonality}
Two complex signals $x(t)$ and $y(t)$ are orthogonal if:
$\int\limits_{-\infty}^{\infty}\!x(t)\cconjugate{y}(t)\,\dt = 0$,
where $\star$ denotes the complex conjugation.
\end{definition}
Orthogonality is invariant to multiplication of either
signal by any real or complex constant.
This is evident from the definition, since if
the integral is zero for a given pair of $x(t)$ and $y(t)$
it must also be zero for $\zeta y(t)$ or $\zeta x(t)$,
where $\zeta$ is a complex constant.
Therefore, given a pair of orthogonal signals,
one of the signals may be rotated by an arbitrary angle
in the complex plane (constant with respect to time),
and still remain orthogonal to the other.
Obviously the same applies for a real constant,
so scaling one or both of the signals does not alter their orthogonality.

Consider the integral in Definition \ref{def:orthogonality}
separated into real and imaginary parts:
\begin{align}
\label{eqn:seporthoreal}
\int\limits_{-\infty}^{\infty}\!(x_r(t)y_r(t)+x_i(t)y_i(t))\,\dt&=0\\
\label{eqn:seporthoimag}
\int\limits_{-\infty}^{\infty}\!(x_i(t)y_r(t)-x_r(t)y_i(t))\,\dt&=0
\end{align}
where $x_r(t)$ (resp. $y_r(t)$) is the real part of $x(t)$
(resp. $y(t)$) and $x_i(t)$ (resp. $y_i(t)$) is the imaginary part of $x(t)$
(resp. $y(t)$). 
There are two significantly different ways in which these integrals
can vanish.
One way is for all four of the products
within them to vanish when integrated:
\begin{align*}
\int\limits_{-\infty}^{\infty}\!x_r(t)y_r(t)\,\dt=0\qquad
\int\limits_{-\infty}^{\infty}\!x_i(t)y_i(t)\,\dt=&0\\
\int\limits_{-\infty}^{\infty}\!x_i(t)y_r(t)\,\dt=0\qquad
\int\limits_{-\infty}^{\infty}\!x_r(t)y_i(t)\,\dt=&0
\end{align*}
This can only happen if both complex signals have an imaginary
part orthogonal to their real part\footnote{An example of such
a signal can be created numerically by taking four singular vectors
from a singular value decomposition of a real matrix, since the
singular vectors are mutually orthogonal.}. In this case only,
orthogonality is invariant to complex conjugation of either or
both of the complex signals.

In more general cases the integrals in \eqref{eqn:seporthoreal}
and \eqref{eqn:seporthoimag} vanish because the product terms
within them integrate to equal values and appropriate signs
(opposite signs in the first case, the same sign in the second).

Of course, a hybrid of these possibilities is also possible,
and indeed occurs in the next section: one integral vanishes
because each of the product terms within it integrates to zero,
and the other vanishes because the product terms integrate to
the same value with appropriate sign, and cancel.

\subsection{Hilbert transform}
\label{sec:hilbert}
\begin{figure}
\setlength{\unitlength}{0.01\columnwidth}
\thicklines
\begin{picture}(100,30)
\put(70,13){\circle{8}}
\put(70,13){\makebox(0,0){$\times$}}
\put(70,4){\vector(0,1){5}}
\put(70,0){\makebox(0,3){$-1$}}
\put(30,20){\framebox(20,10){$\mathcal{H}$}}
\put(30,8){\framebox(20,10){$\mathcal{H}$}}
\put(15,13){\vector(1,0){15}}
\put(15,25){\vector(1,0){15}}
\put(50,13){\vector(1,0){16}}
\put(50,25){\vector(1,0){35}}
\put(74,13){\vector(1,0){10}}
\put(8,25){\makebox(0,0){$z_r(t)$}}
\put(8,13){\makebox(0,0){$z_i(t)$}}
\put(85,25){\makebox(15,0){$o_r(t)$}}
\put(85,13){\makebox(15,0){$o_i(t)$}}
\end{picture}
\caption{\label{fig:dual_hilbert}Construction of the orthogonal signal using a pair of Hilbert transforms.}
\end{figure}
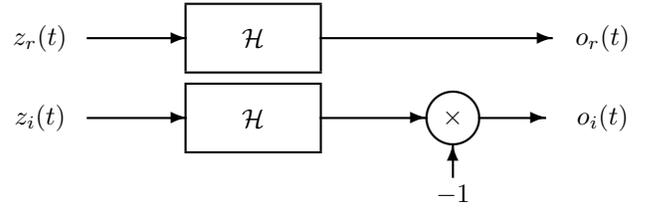

An orthogonal signal may be constructed using a pair of
classical Hilbert transforms as shown in Figure~\ref{fig:dual_hilbert}.
\begin{lemma}
\label{lemma:hilbertcross}
Given two arbitrary real-valued functions $f(t)$ and $g(t)$
the following holds, where $\hilbert{x(t)}$ is the Hilbert transform
of $x(t)$, provided the Hilbert transforms exist:
\begin{equation}
\label{eqn:hilbertcross}
\int\limits_{-\infty}^{\infty}f(t)\hilbert{g(t)}\,\dt
\,+\!
\int\limits_{-\infty}^{\infty}\hilbert{f(t)}g(t)\,\dt                             
=0
\end{equation}
\end{lemma}
\begin{IEEEproof}
We make use of the inner product form of Parseval's Theorem
\cite[\S\,27.9 equation (27.34)]{Jordan:2008} or
\cite[Theorem 7.5]{Champeney:1987}:
\begin{equation}
\label{eqn:innerparseval}
\int\limits_{-\infty}^{\infty}x(t)\cconjugate{y}(t)\,\dt
=
\int\limits_{-\infty}^{\infty}X(\nu)\cconjugate{Y}(\nu)\,\dnu
\end{equation}
where $X(\nu)$ is the Fourier transform of $x(t)$ and similarly for $y(t)$.
Note that, if $y(t)$ is real, as it is in what follows, $Y(\nu)$ must still
be conjugated.
We also make use of the Fourier transform of a Hilbert transform
\cite[(2.13-8, p\,68)]{Stark:1979}:
$\ffourier{\hilbert{x(t)}}=-i\sign(\nu)X(\nu)$, where $X(\nu)=\ffourier{x(t)}$.
Using \eqref{eqn:innerparseval}, we rewrite \eqref{eqn:hilbertcross} as:
\[
\int\limits_{-\infty}^{\infty}\!\cconjugate{F}(\nu)i\sign(\nu)G(\nu)\,\dnu
+
\!\!\int\limits_{-\infty}^{\infty}\!i\sign(\nu)F(\nu)\cconjugate{G}(\nu)\,\dnu=0
\]
which may be simplified to:
\[
i\!\int\limits_{-\infty}^{\infty}\!\sign(\nu)\left[\cconjugate{F}(\nu)G(\nu)
                                                  +F(\nu)\cconjugate{G}(\nu)
                                            \right]\,\dnu=0
\]
Since $\cconjugate{F}(\nu)G(\nu)$ and $F(\nu)\cconjugate{G}(\nu)$ are
conjugates, their imaginary parts cancel and the term inside the brackets
is real.
Since $F(\nu)$ and $G(\nu)$ are the Fourier transforms of real signals,
their real parts are even functions of $\nu$.
Therefore when multiplied by the signum function, they integrate to zero.
\end{IEEEproof}
\begin{proposition}
\label{prop:dual_hilbert}
A complex signal $o(t)$ constructed from an arbitrary complex signal
$z(t)$ using two classical Hilbert transforms operating independently
on the real and imaginary parts,
as represented in Figure~\ref{fig:dual_hilbert} and below:
\begin{equation}
\label{eqn:orthogonal}
o(t) = \hilbert{z_r(t)} - i\hilbert{z_i(t)}
\end{equation}
is orthogonal to $z(t)$, that is:

$\int\limits_{-\infty}^{\infty}z(t)\cconjugate{o}(t)\,\dt = 0$

\end{proposition}

\begin{IEEEproof}
First consider the real part of the integral:
\begin{equation*}
\int\limits_{-\infty}^{\infty}\left[z_r(t)o_r(t)+z_i(t)o_i(t)\right]\,\dt = 0
\end{equation*}
and write it as two integrals with the real and imaginary parts of $o(t)$
replaced by their equivalents from \eqref{eqn:orthogonal}:
\begin{equation*}
\int\limits_{-\infty}^{\infty}z_r(t)\hilbert{z_r(t)}\,\dt\,-\!
\int\limits_{-\infty}^{\infty}z_i(t)\hilbert{z_i(t)}\,\dt =0
\end{equation*}
Both integrals evaluate to zero as a consequence of orthogonality between:
$z_r(t)$ and its Hilbert transform;
and $z_i(t)$ and its Hilbert transform.
Now consider the imaginary part of the integral:
\begin{equation*}
\int\limits_{-\infty}^{\infty}\left[z_i(t)o_r(t)-z_r(t)o_i(t)\right]\,\dt =0
\end{equation*}
and replace the real and imaginary parts of $o(t)$ as before:
\begin{equation*}
\int\limits_{-\infty}^{\infty}z_i(t)\hilbert{z_r(t)}\,\dt\,+\!
\int\limits_{-\infty}^{\infty}z_r(t)\hilbert{z_i(t)}\,\dt =0                           
\end{equation*}
This is true by Lemma~\ref{lemma:hilbertcross}. 

\end{IEEEproof}

\subsection{Complex Fourier transform pair}
An alternative method to construct an orthogonal signal is shown in
Figure~\ref{fig:complex_ft} using a complex Fourier transform pair operating
directly on the complex signal, rather than Hilbert transforms
operating on the real and imaginary parts separately.
The Fourier spectrum of $z(t)$ is multiplied by a signum function times
$i$, equivalent to multiplication of positive frequencies by $e^{i\pi/2}$
and negative frequencies by $e^{-i\pi/2}$.
This has the effect of rotating positive frequency coefficients by
$\pi/2$ and negative frequency coefficients by $-\pi/2$.
An inverse Fourier transform then yields the conjugate of the orthogonal
signal.
\begin{proposition}
\label{prop:dual_fft}
A complex signal $o(t)$ constructed from an arbitrary complex signal
$z(t)$ using a complex Fourier transform as represented in
Figure~\ref{fig:complex_ft} and below, is identical to the orthogonal
signal in Proposition~\ref{prop:dual_hilbert}, equation \eqref{eqn:orthogonal}.
\begin{equation}
\label{eqn:dual_fft_orthogonal}
\cconjugate{o}(t) = \ifourier{-i\sign(\nu)Z(\nu)}
\end{equation}
where $\fourier$ denotes a Fourier transform and $Z(\nu)=\ffourier{z(t)}$.
\end{proposition}
\begin{IEEEproof}
From equation \eqref{eqn:orthogonal}, $\cconjugate{o}(t)=\hilbert{z_r(t)} + i\hilbert{z_i(t)}$.
Taking the Fourier transforms of the Hilbert transforms (see Lemma \ref{lemma:hilbertcross}),
we have:
\begin{align*}
\ffourier{\cconjugate{o}(t)}&=-i\sign(\nu)Z_r(\nu)-i\sign(\nu)Z_i(\nu)\\
                            &=-i\sign(\nu)Z(\nu)
\end{align*}
from the linearity of the Fourier transform.

\end{IEEEproof}
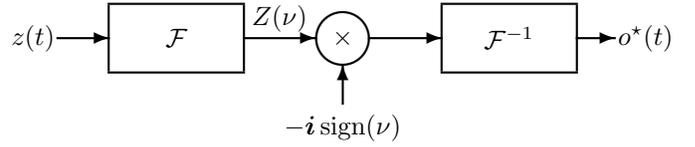
\begin{figure}
\setlength{\unitlength}{0.01\columnwidth}
\thicklines
\begin{picture}(100,21)
\put(50,16){\circle{8}}
\put(50,16){\makebox(0,0){$\times$}}
\put(50,6){\vector(0,1){6}}
\put(50,0){\makebox(0,5){$-\i\sign(\nu)$}}
\put(15,11){\framebox(20,10){$\fourier$}}
\put(65,11){\framebox(20,10){$\fourier^{-1}$}}
\put(0,16){\makebox(7,0){$z(t)$}}
\put(7,16){\vector(1,0){8}}
\put(35,16){\vector(1,0){11}}
\put(35,16){\makebox(11,6){$Z(\nu)$}}
\put(54,16){\vector(1,0){11}}
\put(85,16){\vector(1,0){6}}
\put(91,16){\makebox(9,0){$\cconjugate{o}(t)$}}
\end{picture}
\caption{\label{fig:complex_ft}Construction of the orthogonal signal using a complex Fourier transform.}
\end{figure}

\subsection{Discussion}
We have shown in this section how to construct a complex signal $o(t)$
orthogonal to an arbitrary complex signal $z(t)$ using classical Hilbert
or complex Fourier transforms.
The notion of analytic signal, however, requires more than just an
orthogonal signal --- generalising from the classical case, we must
treat the original and orthogonal complex signals as a pair in order
to define the generalisations of envelope and phase.
In the next section we introduce the quaternion Fourier transform,
and show that this permits us to construct a \emph{hyperanalytic signal}
from which both the original signal and the orthogonal signal
defined in the current section may be extracted.
Using the quaternion hyperanalytic signal we then show how to
construct the envelope and phase in the complex case.

\section{Construction of an orthogonal signal using a quaternion Fourier Transform}
\label{QFTsection}

In this section, we will be concerned with the definition and
properties of the quaternion Fourier transform (QFT) of complex valued
signals. Before introducing the main definitions, we give some
prerequisites on quaternion valued signals.

\subsection{Preliminary remarks}
Quaternions were discovered by Sir W.R. Hamilton in 1843
\cite{Hamilton:1853}. Quaternions are 4D hypercomplex numbers that form a noncommutative
division ring denoted $\H$. A quaternion $q \in \H$
has a Cartesian form: $q=a+b\i+c\j+d\k$, with $a,b,c,d
\in {\mathbb R}$ and $\i,\j,\k$ roots of $-1$ satisfying
$\i^2=\j^2=\k^2=\i\j\k=-1$.
The {\em scalar part} of $q$ is $a$: ${\cal S}(q)=a$.
The {\em vector part} of $q$ is: ${\cal V}(q)=q-{\cal S}(q)$.
Quaternion multiplication is not commutative, so that in
general $qp\neq pq$ for $p,q \in \H$.
The conjugate of $q$ is $\qconjugate{q}=a-b\i-c\j-d\k$.
The norm of $q$ is $\norm{q}=|q|^2=(a^2+b^2+c^2+d^2)=q\qconjugate{q}$.
A quaternion with $a=0$ is called pure.
If $|q|=1$, then $q$ is called a {\em unit} quaternion.
The inverse of $q$ is $q^{-1}=\qconjugate{q}/\norm{q}$.
Pure unit quaternions are special quaternions,
among which are $\i$, $\j$ and $\k$.
Together with the identity of $\H$, they form a {\em quaternion basis}:
$(1,\i,\j,\k)$. In fact, given any two unit pure quaternions, $\mu$
and $\xi$, which are orthogonal to each other ({\em i.e.} ${\cal
  S}(\mu\xi)=0$), then $(1,\mu,\xi,\mu\xi)$ is a quaternion basis.

Quaternions can also be viewed as complex numbers with complex
components, {\em i.e.} one can write $q=z_1+\j z_2 $ in the
basis $(1,\i,\j,\k)$ with $z_1,z_2 \in \C^{\i}$, {\em i.e.} $z_{\alpha}=\Re(z_{\alpha})+\i
\Im(z_{\alpha})$ for $\alpha=1,2$\footnote{In the sequel, we note as $\C^{\mu}$
a complex number with root of $-1$ being $\mu$. Note that
these are degenerate quaternions.}. This is called the
Cayley-Dickson form. Now, it is possible to define a
generalized Cayley-Dickson form of a given quaternion using an
arbitrary quaternion basis $(1,\mu,\xi,\mu\xi)$. A quaternion $q$ can
then be written as $q=w_1+ \xi w_2 $ where $w_1,w_2 \in \C^{\mu}$.
  
Any quaternion $q$ also has a unique polar Cayley-Dickson
  form \cite{10.1007/s00006-008-0128-1} given by: $q=A \exp(B\j)$, where $A$ is the complex modulus of $q$
  and $B$ is its complex phase. Expressions of $A$ and $B$ can be
  found in \cite{10.1007/s00006-008-0128-1}. This notation will be of
  use in Section \ref{PhaseAttribsection}.

\subsection{1D Quaternion Fourier transform}
In this paper, we use a 1D version of the (right) QFT first
defined in discrete-time form in \cite{SangwineEll:1998}. 
Thus, it is necessary to specify the axis (a pure unit
quaternion) of the transform. So, we will denote $\qfourier{\mu}{\,}$ a QFT with
transformation axis $\mu$. This will become clear from the definition
below. We now present the definition and some properties of the
transform used here.

\begin{definition}
\label{def:qft}
Given a complex valued signal $z(t)=z_r(t)+\xi z_i(t)$, its quaternion
Fourier transform with respect to axis $\mu$ is:
\begin{equation}
Z_{\mu}(\nu)=\qfourier{\mu}{z(t)}=\int\limits_{-\infty}^{+\infty}z(t)e^{-\mu 2 \pi \nu t} \dt
\label{Eqn:QFT2}
\end{equation}
and the inverse transform is:
\begin{equation}
z(t)=\iqfourier{\mu}{Z_{\mu}(\nu)}=\int\limits_{-\infty}^{+\infty}Z_{\mu}(\nu)e^{\mu 2 \pi \nu t} d\nu
\label{Eqn:QFT3}
\end{equation}
where in both cases,
the transform axis $\mu$ can be chosen to be any unit pure quaternion orthogonal to $\xi$
so that $\left(1,\xi,\mu,\xi\mu\right)$ is a quaternion basis.
\end{definition}

\begin{property}
\label{prop:symQFT}
Given a complex signal $z(t)=z_r(t)+\xi z_i(t)$ and its quaternion
Fourier transform denoted $Z_{\mu}(\nu)$, then the following properties hold:
\begin{itemize}
\item The even part of $z_r(t)$ is in          $\Re(Z_{\mu}(\nu))$.
\item The odd  part of $z_r(t)$ is in    $\Im_{\mu}(Z_{\mu}(\nu))$
\item The even part of $z_i(t)$ is in    $\Im_{\xi}(Z_{\mu}(\nu))$.
\item The odd  part of $z_i(t)$ is in $\Im_{\xi\mu}(Z_{\mu}(\nu))$
\end{itemize}
\end{property}
\begin{IEEEproof}
Expand \eqref{Eqn:QFT2} into real and imaginary parts with respect to $\xi$,
and expand the quaternion exponential into cosine and sine components:
\begin{align*}
Z_{\mu}(\nu)
&=\int\limits_{-\infty}^{+\infty}{
                                  \!\left[z_r(t)+\xi z_i(t)\right]
                                  }
                                   \left[\cos(2\pi\nu t) - \mu\sin(2\pi\nu t)\right]\dt\\ 
&=\phantom{+\xi}\!\int\limits_{-\infty}^{+\infty}\!z_r(t)\cos(2\pi\nu t)\dt 
            -\mu\!\int\limits_{-\infty}^{+\infty}\!z_r(t)\sin(2\pi\nu t)\dt\\
&\phantom{=}+\xi\!\int\limits_{-\infty}^{+\infty}\!z_i(t)\cos(2\pi\nu t)\dt 
         -\xi\mu\!\int\limits_{-\infty}^{+\infty}\!z_i(t)\sin(2\pi\nu t)\dt 
\end{align*}
from which the stated properties are evident.
\end{IEEEproof}

These properties are central to the justification of the use of the QFT to
analyze a complex valued signal carrying complementary but different
information in its real and imaginary parts. Using the QFT, it is
possible to have the odd and even parts of the real and imaginary parts of
the signal in four different components in the transform
domain. This idea was also the initial motivation of B\"ulow,
Sommer and Felsberg when they developed the monogenic signal for images
\cite{Bulow:1999,BulowSommer:2001,FelsbergSommer:2001}. Note that the
use of hypercomplex Fourier transforms was
originally introduced in 2D Nuclear Magnetic Resonance image analysis \cite{Ernst:1987,Delsuc:1988}.

We now turn to the link between a complex signal and the quaternion
signal that can be uniquely associated to it.
\begin{property}
A one-sided QFT $X(\nu)$ with $X(\nu)=0$ ({\em i;e.}$\forall \nu < 0$)
corresponds to a quaternion valued signal in the time domain. 
\end{property}
\begin{IEEEproof}
A complex signal with values in $\C^{\xi}$ has a $QFT_{\mu}$ with the symmetry
properties listed in property \ref{prop:symQFT}.
Now, cancelling its negative
frequencies is equivalent to adding a signal with a
QFT with the following
properties: $\Re$ part is odd, $\Im_{\mu}$ is odd, $\Im_{\xi}$ is even
and $\Im_{\xi\mu}$ is even. By linearity of the inverse QFT,
the former two parts of the QFT produce a complex signal in the time
domain with $\Re$ and $\Im_{\xi}$ parts, while the two latter produce a
quaternion signal in the time domain with $\Im_{\mu}$ and $\Im_{\xi\mu}$
parts. Together, and again by linearity, the resulting signal in the
time domain is a full 4D quaternion signal.
\end{IEEEproof}

\begin{property}
Given a complex signal $x(t)$, one can associate to it a unique
canonical pair corresponding to a modulus and phase. These
modulus and phase are uniquely defined through the hyperanalytic
signal, which is quaternion valued.
\end{property}
\begin{IEEEproof}
Cancelling the negative frequencies of the QFT leads to a quaternion
signal in the time domain. Then, any quaternion signal has a modulus
and phase defined using its CD polar form.
\end{IEEEproof}

\subsection{Convolution}
We consider the special case of convolution of a complex signal by a
real signal. Consider $g$ and $f$ such that: $g:{\mathbb
R}^+\rightarrow{\mathbb C}^{\xi}$ and $f:{\mathbb R}^+\rightarrow{\mathbb
R}$. Now, consider the $\QFT_{\mu}$ of their convolution:
\begin{equation}
\begin{aligned}
\qfourier{\mu}{g*f(t)} &= \int\limits_{-\infty}^{+\infty}
                     \int\limits_{-\infty}^{+\infty}
                     g(\tau)f(t-\tau)\dtau
                     e^{-2\mu \pi \nu t} \dt \\
                  &= \int\limits_{-\infty}^{+\infty}
                     \int\limits_{-\infty}^{+\infty}
                     g(\tau)e^{-\mu 2 \pi \nu (t'+\tau)}f(t')d\tau \dt'\\
                  &= \int\limits_{-\infty}^{+\infty}
                     g(\tau)e^{-2\mu \pi \nu \tau}\dtau
                     \int\limits_{-\infty}^{+\infty}
                     f(t')e^{-2\mu \pi \nu t'}\dt'\\
                  &= \qfourier{\mu}{g(t)}\qfourier{\mu}{f(t)}\\
                  &= \qfourier{\mu}{f(t)}\qfourier{\mu}{g(t)}
\end{aligned}
\label{Eqn:QFT3A}
\end{equation}
Thus, the definition used for the QFT here verifies the convolution
theorem in the considered case. This specific case will be of use in
our definition of the hyperanalytic signal. 

\subsection{The quaternion Fourier transform of the Hilbert transform}

It is straightforward to verify that the quaternion Fourier transform
of a real signal $x(t)=\frac{1}{\pi t}$ is $-\mu\sign(\nu)$,
where $\mu$ is the axis of the transform.
Substituting $x(t)$ into \eqref{Eqn:QFT2}, we get:
\begin{equation*}
\qfourier{\mu}{\frac{1}{\pi t}}=\frac{1}{\pi}\int\limits_{-\infty}^{+\infty}\frac{e^{-\mu 2\pi\nu t}}{t}\,\dt
\end{equation*}
and this is clearly isomorphic to the classical complex case. The solution
in the classical case is $-i\sign(\nu)$, and hence in the quaternion case
must be as stated above.

It is also straightforward to see that, given an arbitrary real signal $y(t)$,
subject only to the constraint that its classical Hilbert transform
$\hilbert{y(t)}$ exists, then one can easily show that the classical Hilbert transform
of the signal may be obtained using a quaternion Fourier transform as
follows:
\begin{equation}
\label{eqn:hilbertbyqft}
\hilbert{y(t)} = \iqfourier{\mu}{-\mu\sign(\nu)Y_{\mu}(\nu)}
\end{equation}
where $Y_{\mu}(\nu)=\qfourier{\mu}{y(t)}$.
This result follows from the isomorphism between the quaternion and complex
Fourier transforms when operating on a real signal, and it may be seen to be
the result of a convolution between the signal $y(t)$ and the quaternion
Fourier transform of $x(t)=\frac{1}{\pi t}$.
Note that $\mu$ and $Y_{\mu}(\nu)$ commute as a consequence of $y(t)$ being
real.

\section{The hyperanalytic signal as a quaternion signal with a one-sided spectrum}
\label{qonesidedsection}

We define the hyperanalytic signal $h(t)$ by a similar approach
to that originally developed by Ville \cite{Ville:1948}. 
The following definitions give the details of the construction of this signal.
\emph{Note that the signal $z(t)$ is considered to be non-analytic in the classical (complex) sense,
that is its real and imaginary parts are \emph{not}
orthogonal. However, the following definitions are valid if $z(t)$ is
analytic, as it can be considered as a degenate case of the more
general non-analytic case.} 

\begin{definition}
\label{def:hyperhilbert}
Consider a complex signal $z(t)=z_r(t)+\xi z_i(t)$
and its quaternion Fourier transform $Z_{\mu}(\nu)$ as defined in Definition~\ref{def:qft}.
Then, the hypercomplex analogue of the Hilbert transform of $z(t)$, is as follows:
\begin{equation}
\qhilbert{\mu}{z(t)} =\iqfourier{\mu}{-\mu\sign(\nu)Z_{\mu}}(\nu)
\end{equation}
where the Hilbert transform is thought of as:
$\qhilbert{\mu}{z(t)}=p.v.\left(z*\frac{1}{\pi t}\right)$, where the
principal value (p.v.) is understood in its classical way. This result follows from equation \eqref{eqn:hilbertbyqft} and the linearity
of the quaternion Fourier transform.
Notice that the result will take the same form as $z(t)$,
namely a quaternion signal $o_r(t) + \xi o_i(t)$
isomorphic to the complex signal $\cconjugate{o}(t)=o_r(t) +   i o_i(t)$.
To extract the imaginary part, the vector part of the quaternion signal must
be multiplied by $-\xi$.
An alternative is to take the scalar or inner product of the vector part with
$\xi$.
Note that $\mu$ and $Z_{\mu}(\nu)$ anticommute because $\mu$ is orthogonal
to $\xi$, the axis of $Z_{\mu}(\nu)$.
Therefore the ordering is not arbitrary, but changing it simply conjugates
the result.
\end{definition}

\begin{definition}
\label{def:hyperanalytic}
Given a complex valued signal $z(t)$ that can be expressed in the form of a quaternion as $z(t)=z_r(t)+\xi z_i(t)$, and given a pure unit quaternion $\mu$ such that $(1,\xi,\mu,\xi\mu)$ is a quaternion basis, then the hyperanalytic signal of $z(t)$, denoted $h(t)$ is given by: 
\begin{equation}
h(t)=z(t)+\mu\qhilbert{\mu}{z(t)}
\label{defHana}
\end{equation}
where $\qhilbert{\mu}{z(t)}$ is the hypercomplex analogue of the
Hilbert transform of $z(t)$ defined in the preceding definition.
The quaternion Fourier transform of this hyperanalytic signal is thus:
\begin{align*}
H_{\mu}(\nu)
&=Z_{\mu}(\nu)-\mu^{2}\sign(\nu)Z_{\mu}(\nu)\\
&=\left[1+\sign(\nu)\right]Z_{\mu}
\end{align*}
which is a direct extension of the `classical' analytic signal.
\end{definition}

This result is unique to the
quaternion Fourier transform representation of the hyperanalytic signal --- the
hyperanalytic signal has a one-sided quaternion Fourier spectrum.
This means that the hyperanalytic signal may be constructed from a complex signal
$z(t)$ in exactly the same way that an analytic signal may be constructed from a
real signal $x(t)$, by suppression of negative frequencies in the Fourier spectrum.
The only difference is that in the hyperanalytic case,
a quaternion rather than a complex Fourier transform must be used,
and of course the complex signal must be put in the form $z(t)=z_r(t)+\xi z_i(t)$
which, although a quaternion signal, is isomorphic to the original complex signal.

A second important property of the hyperanalytic signal is that it
maintains a separation between the different even and odd parts of
the original signal.

\begin{property}
The original signal $z(t)$ is the \emph{simplex} part \cite[Theorem 1]{10.1109/ICIP.2000.899828},
\cite[\S\,13.1.3]{Multivariate:2009} of its corresponding hyperanalytic signal. The perplex part is the orthogonal or `quadrature' component, $o(t)$. 
They are obtained by:
\begin{align}
\label{eqn:simplex}
z(t)&=\frac{1}{2}\left(h(t) - \xi h(t)\xi\right)\\
\label{eqn:perplex}
o(t)&=\frac{1}{2}\left(h(t) + \xi h(t)\xi\right)
\end{align}
\end{property}
\begin{IEEEproof}
This follows from equation \eqref{defHana}.
Writing this in full by substituting the orthogonal signal for $\qhilbert{\mu}{z(t)}$:
\[
h(t) = z(t) + \mu o(t) = z_r(t) + \xi z_i(t) + \mu o_r(t) + \mu\xi o_i(t) 
\]
and substituting this into equation \eqref{eqn:simplex}, we get:
\begin{align*}
z(t)&=\frac{1}{2}\left(
\begin{aligned}&\phantom{\xi\left[\right.}z_r(t) + \xi z_i(t) + \mu o_r(t) + \mu\xi o_i(t)\\
              -&         \xi\left[        z_r(t) + \xi z_i(t) + \mu o_r(t) + \mu\xi o_i(t) \right]\xi
\end{aligned}\right)\\
\intertext{and since $\xi$ and $\mu$ are orthogonal unit pure quaternions, $\xi\mu=-\mu\xi$:} 
    &=\frac{1}{2}\left(
\begin{aligned}&\phantom{\xi\left[\right.}z_r(t) + \xi z_i(t) + \mu o_r(t) + \mu\xi o_i(t)\\
              +&\phantom{\xi\left[\right.}z_r(t) + \xi z_i(t) - \mu o_r(t) - \mu\xi o_i(t)
\end{aligned}\right)
\end{align*}
from which the first part of the result follows. Equation \eqref{eqn:perplex} differs
only in the sign of the second term, and it is straightforward to see that if $h(t)$
is substituted, $z(t)$ cancels out, leaving $o(t)$. 
\end{IEEEproof}

\section{Phase and envelope from quaternion representation}
\label{PhaseAttribsection}

As was shown in §\,\ref{sec:hilbert} and Figure~\ref{fig:dual_hilbert},
it is possible to construct an orthogonal signal using classical
Hilbert transforms operating independently on the real and imaginary
parts of a complex signal $z(t)$.
However, to extend the analytic signal concept to the hyperanalytic
case, we need a method to construct the complex envelope
(also known as the \emph{instantaneous amplitude})
and the complex phase
(also known as the \emph{instantaneous phase})
from a \emph{pair} of complex signals
(the original complex signal $z(t)$,
and the complex orthogonal signal $o(t)$).

A complex envelope may be obtained,
\emph{apart from the signs of the real and imaginary parts,
and therefore the quadrant},
from the moduli of the analytic signals
of the real and imaginary parts of the original signal. The resulting envelope will always be in the first quadrant and
therefore in general it will not be a correct instantaneous amplitude.

In this paper the method we present is based on quaternion algebra
which is a natural choice given that a quaternion may be expressed
as a pair of complex numbers in the so-called Cayley-Dickson form,
and, as we have shown earlier, the hyperanalytic signal may be
computed using a discrete quaternion Fourier transform.

\subsection{Complex envelope}

In a previous paper \cite{SangwineLeBihan:2007} we showed how the envelope
of the hyperanalytic signal could be obtained from a biquaternion Fourier
transform, using simply the complex modulus of the biquaternion
hyperanalytic signal.
This approach has a limitation, namely that the envelope
is always in the right half of the complex plane.
(This is because the modulus of the biquaternion hyperanalytic signal is
obtained as the square root of the sum of the squares of the components,
and the complex square root is conventionally in the right half-plane.)
A second, more serious, limitation of the biquaternion Fourier transform is that under
some circumstances the complex envelope can vanish because the biquaternion
hyperanalytic signal is \emph{nilpotent} \cite{10.1007/s00006-010-0202-3}:
that is it has values with zero modulus,
even though the biquaternion values are non-zero.

The mathematical key to constructing the complex envelope from a quaternion hyperanalytic
signal is given in a paper written in 2008 \cite{10.1007/s00006-008-0128-1}.
In that paper we set out a new polar representation for quaternions that has a
complex modulus and a complex argument:
\begin{equation}
\label{eqn:cdpolar}
q = A \exp\left(B\j\right)
\end{equation}
where $A$ and $B$ are complex
(more accurately they are degenerate quaternions of the form $w + x\i + 0\j + 0\k$).

The representation as published is based on the standard quaternion basis
$(1,\i,\j,\k)$, although clearly it could use any other basis.
For simplicity of numerical computation, we use the standard basis.

This representation of a quaternion is not trivial to compute numerically
and one of the main problems we had in arriving at the results presented
in this paper was to refine the computation of the representation in
\eqref{eqn:cdpolar} to minimise numerical problems.
We do not discuss these issues here, as they are solely concerned with
correct numerical implementation of \eqref{eqn:cdpolar}.
The code is open source, and the numerical issues are documented in the code
\cite[Function: \texttt{cdpolar.m}]{qtfm}.
It is also necessary to use phase unwrapping on $B$ to obtain a continuous
envelope and phase without discontinuities when computing $A(t)$ and $B(t)$
using \eqref{eqn:cdpolar} from a hyperanalytic signal $h(t)$, due to
an inherent ambiguity of sign in $A$ and $B$ \cite[§2.1]{10.1007/s00006-008-0128-1}.
(In a similar way the complex envelope obtained from a biquaternion hyperanalytic
signal may be phase unwrapped to obtain an envelope in all four quadrants.)

\begin{definition}
\label{def:cdenv}
Given a complex signal $z(t)$, and its hyperanalytic signal $h(t)$
constructed according to Definition~\ref{def:hyperanalytic},
expressed in the standard quaternion basis $(1,\i,\j,\k)$,
the instantaneous amplitude, or complex envelope, $e(t)$, and the
instantaneous phase, $\phi(t)$, are given by the polar form of a
quaternion defined in \eqref{eqn:cdpolar}:
\[
h(t) = z(t) + o(t)\j = e(t)\exp(\phi(t)\j)
\]
where $z(t)$, $o(t)$, $e(t)$, and $\phi(t)$ are complex.
\end{definition}

\subsection{Complex modulation of a real carrier}
\label{Modulationsection}

\begin{figure}
\includegraphics[width=\columnwidth]{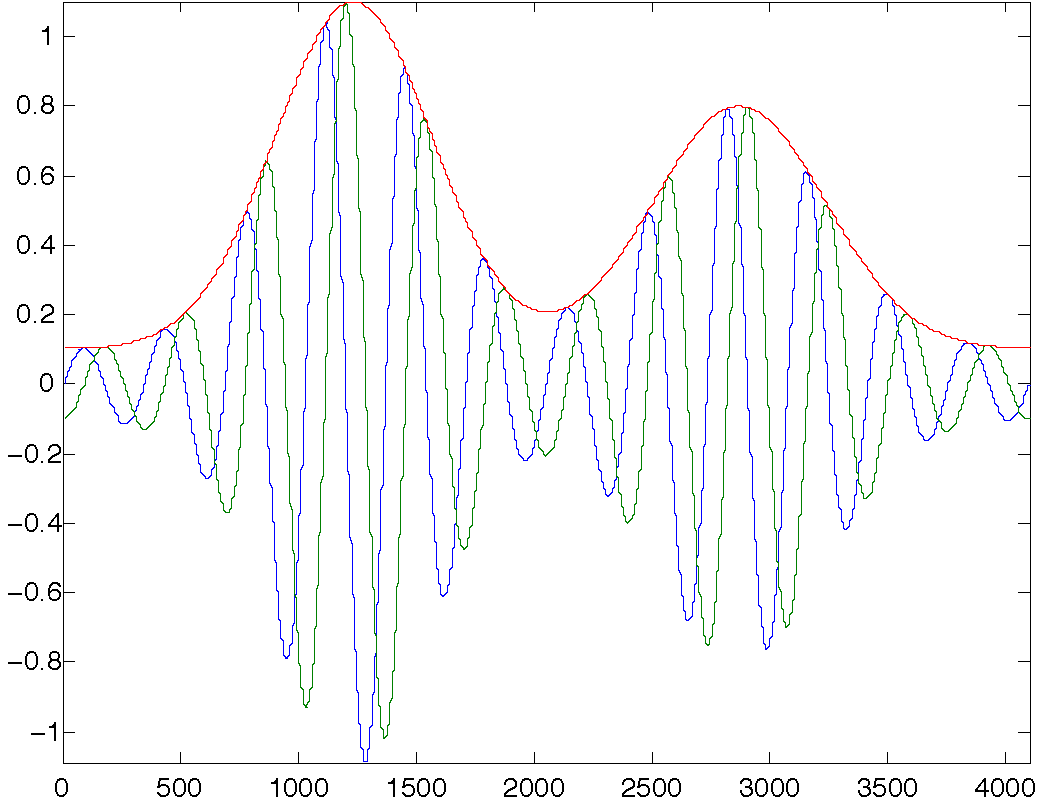}
\caption{\label{fig:analytic}Classical envelope of an amplitude modulated signal obtained as the modulus of the analytic signal. Blue: original modulated signal; green: orthogonal signal; red: envelope.}
\end{figure}

A classic result in the case of the analytic signal concerns
amplitude modulation, where the envelope of the analytic
signal (that is its modulus) allows the modulating signal
to be recovered, subject to certain constraints such as the
depth of modulation \cite[\S\,7.16]{Hahn:Poularikas}.
Figure \ref{fig:analytic} shows an example.
In this section we extend this result to the complex case.
\begin{property}
\label{prop:instphase}
Given a signal $z(t)$ consisting of a real sinusoidal `carrier' $f(t) = A\cos(2\pi\omega t+\theta)$,
and a \emph{complex} `modulating' signal $g(t)$:
\[
z(t) = g(t)\,f(t)
\]
where the frequency content of $g(t)$ is below the frequency of the carrier, $\omega$,
then the instantaneous phase $\phi(t)$ as defined in Definition \ref{def:cdenv}
is real and equal to the instantaneous phase of $f(t)$, that is $2\pi\omega t + \theta$,
and the envelope $e(t)$ is a scaled version of the complex modulating signal, $A\,g(t)$.
\end{property}
\begin{IEEEproof}
Separate $z(t)$ into real and imaginary parts:
\begin{align*}
z_r(t) &= g_r(t)A\cos(2\pi\omega t+\theta)\\
z_i(t) &= g_i(t)A\cos(2\pi\omega t+\theta)
\end{align*}
By Bedrosian's Theorem \cite{Bedrosian:1963}\cite[\S\,7.12 and Table 7.7.1]{Hahn:Poularikas},
the Hilbert transforms of the real and imaginary parts of $z(t)$ are:
\begin{align*}
\hilbert{z_r(t)} &= g_r(t)A\sin(2\pi\omega t+\theta)\\
\hilbert{z_i(t)} &= g_i(t)A\sin(2\pi\omega t+\theta)
\end{align*}
We can therefore write the hyperanalytic signal as follows:
\begin{align*}
h(t) &= z_r(t) + \i z_i(t) + \j\hilbert{z_r(t)} + \k\hilbert{z_i(t)}\\
     &= g(t)A[\cos(2\pi\omega t + \theta)
          +\j\sin(2\pi\omega t + \theta)]\\
     &= g(t)A\exp(\j(2\pi\omega t + \theta))
\end{align*}
and comparison with Definition \ref{def:cdenv} shows that $e(t)=A\,g(t)$
and $\phi(t)=2\pi\omega t + \theta$.
\end{IEEEproof}

Note that we have considered above only the case of a real sinusoidal
carrier modulated by a complex modulating signal, and that scaling of the carrier cannot be separated
from scaling of the modulating signal.
A more general case would be a complex carrier, for example,
a complex exponential, or a single sinusoid rotated in the complex
plane, modulated by a real or complex modulating signal, but we now
show that this is ambiguous -- the complex nature of the carrier
cannot be separated from the complex nature of the modulating signal.
Consider a real sinusoidal carrier $f(t)$ modulated by a complex
modulating signal $g(t)$, yielding a modulated signal $m(t)$:
\[
m(t) = f(t)\,g(t)
\]
Now rotate the carrier by an angle $\psi$ in the complex plane so that
it oscillates along a line in the complex plane at angle $\psi$ to the
real axis, instead of along the real axis itself:
\[
m(t) = [f(t)\exp(i\psi)]\,g(t)
\]
Clearly, this modulated signal is indistinguishable from the following,
in which the complex modulating signal $g(t)$ is instead
rotated through the angle $\psi$:
\[
m(t) = f(t)\,[g(t)\exp(i\psi)]
\]
We conclude that the hyperanalytic signal cannot distinguish between a
real carrier modulated by a complex signal and a
complex carrier modulated by a (different) complex
signal: it suffices to consider the former case,
and assume that the modulating signal only is complex.

\section{Properness}

In this section, we study the properness properties of the
hyperanalytic signal. We make use of the term proper here for
deterministic signals, by extension of the concept used for quaternion
random functions in \cite{Vakhania:1998}. The term proper is also
widely used for complex random signals when the real and imaginary parts are orthogonal and have the same
amplitude. Our definition of properness is based on \cite{Amblard:2004} . We use the
  properties of the correlation matrix in the real (Cartesian)
  representation of quaternions. As explained in \cite{Amblard:2004}, a quaternion
can be represented as a four dimensional vector
$q=\left[a,b,c,d\right]^T$. Properness can be observed through the
pattern of the covariance matrix $\Gamma_q$ that contains the
cross-correlations between the four components of the
quaternion. Here, we are interested in the covariance (correlation at zero-lag) of the
hyperanalytic signal $h(t)$. Assuming that the original signal $z(t)$
takes values in $\C^{\xi}$, then $h(t)$ (obtained using a quaternion Fourier transform with axis $\mu$) can be written as:
$$
h(t)=z(t)+o(t)\mu=z_r(t)+ \xi z_i(t) + \left( o_r(t)+\xi o_i(t)\right)\mu
$$
The covariance matrix of $h(t)$, denoted $\Gamma_{h}$, is
thus given as:
$$
\Gamma_{h}=
\begin{bmatrix}
C_{z_r z_r} & C_{z_r z_i} & C_{z_r o_r} & C_{z_r o_i}\\
C_{z_i z_r} & C_{z_i z_i} & C_{z_i o_r} & C_{z_i o_i}\\ 
C_{o_r z_r} & C_{o_r z_i} & C_{o_r o_r} & C_{o_r o_i}\\
C_{o_i z_r} & C_{o_i z_i} & C_{o_i o_r} & C_{o_i o_i}\\
\end{bmatrix}
$$
Now, given a complex signal $z(t)$ taking values in $\C^{\xi}$ and its
hyperanalytic (quaternion valued) signal $h(t)$, then the following stands.
\begin{property}
\label{prop:properC}
If $z(t)$ with values in $\C^{\xi}$ is an improper complex signal, then its
hyperanalytic signal $h(t)$, computed using a quaternion Fourier
transform with axis $\mu$ is $\C^{\mu}$-proper, providing that $(1,\xi,\mu,\xi\mu)$ is a
quaternion basis. As a consequence, the covariance of $h(t)$ is: 
$$
\Gamma_{h}=
\begin{bmatrix}
\alpha & \gamma & \phantom{-}     0 & \phantom{-}\beta \\
\gamma & \omega & \phantom{-}\beta  & \phantom{-}0 \\ 
     0 & \beta  & \phantom{-}\alpha & -\gamma \\
\beta  &     0  &           -\gamma & \phantom{-}\omega\\
\end{bmatrix}
$$
with:
\begin{align*}
\alpha &= \!\int\limits_{-\infty}^{\infty}z^2_r(t)    \dt= \int\limits_{-\infty}^{\infty}o^2_r(t)\dt\\
\omega &= \!\int\limits_{-\infty}^{\infty}z^2_i(t)    \dt= \int\limits_{-\infty}^{\infty}o^2_i(t)\dt\\
\beta  &= \!\int\limits_{-\infty}^{\infty}z_i(t)o_r(t)\dt= \int\limits_{-\infty}^{\infty}z_r(t)o_i(t)\dt\\
\gamma &= \!\int\limits_{-\infty}^{\infty}z_r(t)z_i(t)\dt=-\int\limits_{-\infty}^{\infty}o_r(t)o_i(t)\dt
\end{align*}
\end{property}
\begin{IEEEproof}
We need to prove the zeroing of the terms of $\Gamma_h$ as well as
equalities between some of its terms. Obviously $\Gamma_h$ is a real symmetric
matrix, so we only need proofs for the diagonal and upper-diagonal elements.
We start with the diagonal elements. 

First, recall that, due to Parseval's theorem:
$$
\int_{-\infty}^{\infty}z_r(t)^2\dt=\int_{-\infty}^{\infty}|Z_r(\nu)|^2\dnu
$$
Now by definition of the orthogonal signal $o(t)$, the Fourier
transform of its real part is:
$$
O_r(\nu)=\ffourier{o_r(t)}=\ffourier{\hilbert{z_r(t)}}=-\xi\sign(\nu)Z_r(\nu)
$$
It follows directly that $|O_r(\nu)|^2=|Z_r(\nu)|^2$ and using Parseval's theorem
(which is also valid for $o_r(t)$), the following is true:
$$
\alpha=\int_{-\infty}^{\infty}z^2_r(t)\dt
      =\int_{-\infty}^{\infty}o_r^2(t)\dt
      =\int_{-\infty}^{\infty}\hilbert{z_r(t)}^2\dt
$$
A similar reasoning is valid to show that: 
$$\omega=\int_{-\infty}^{\infty}z_i(t)^2dt=\int_{-\infty}^{\infty}o_i^2(t)\dt$$
 with the small difference that $z_i(t)=-\hilbert{z_i(t)}$, but the sign has no effect on the
reasoning due to the modulus.

We now turn to off-diagonal elements. As $z(t)$ is an {\it improper} signal, its real and imaginary
parts are not orthogonal to each other, so the following is true:
$$
\gamma = \int\limits_{-\infty}^{\infty}z_r(t)z_i(t)dt \neq 0
$$
Now, Parseval's theorem allows us to write:
\[
\int\limits_{-\infty}^{\infty}z_r(t)\cconjugate{\left(\cconjugate{z_i}(t)\right)}\dt=
\int\limits_{-\infty}^{\infty}Z_r(\nu)Z_i(-\nu)\dnu
\]
where $Z_r(\nu)$ and $Z_i(\nu)$ are the classical Fourier transforms
of $z_r(t)$ and $z_i(t)$. We also have that: 
\begin{align*}
o_r(t)&=\,\phantom{-}\hilbert{z_r(t)}\\
o_i(t)&=           - \hilbert{z_i(t)}
\end{align*} 
which implies that:
\begin{alignat*}{2}
O_r(\nu)&=         - \xi\sign(\nu)Z_r(\nu)\\
O_i(\nu)&=\phantom{-}\xi\sign(\nu)Z_i(\nu)
\end{alignat*} 
Therefore, we can write the following equality:
\begin{align*}
\phantom{=}&\int\limits_{-\infty}^{\infty}O_r(\nu)O_i(-\nu)\dnu\\
=&\int\limits_{-\infty}^{\infty}-\xi\sign(\nu)Z_r(\nu)\xi\sign(-\nu)Z_i(-\nu)\dnu\\
=&-\!\int\limits_{-\infty}^{\infty}Z_r(\nu)Z_i(-\nu)\dnu
=-\!\int\limits_{-\infty}^{\infty}z_r(t)z_i(t)\dt
\end{align*}
This shows, due to Parseval's theorem, that: 
$$
\int\limits_{-\infty}^{\infty}o_r(t)o_i(t)\dt=-\!\int\limits_{-\infty}^{\infty}z_r(t)z_i(t)\dt=-\gamma
$$
Then, we have already shown in the proof of Proposition
\ref{prop:dual_hilbert} that:
\begin{equation*}
\int\limits_{-\infty}^{\infty}z_r(t)o_r(t)\dt = 0
\quad\text{and}\quad
\int\limits_{-\infty}^{\infty}z_i(t)o_i(t)\dt = 0
\end{equation*}
due to orthogonality between each signal and its Hilbert transform. 

Now, it was also shown in Proposition \ref{prop:dual_hilbert} that, due to orthogonality between $z(t)$ and $o(t)$ the following stands:
$$\int\limits_{-\infty}^{\infty}z_i(t)o_r(t)\dt=\int\limits_{-\infty}^{\infty}z_r(t)o_i(t)\dt=\beta$$
We now show that in the case where $z(t)$ is {\it improper}, then $\beta \neq 0$.
We start by recalling that:
\begin{align*}
 \int\limits_{-\infty}^{\infty}z_i(t)o_r(t)\dt &
=\int\limits_{-\infty}^{\infty}Z_i(\nu)\cconjugate{O_r}(\nu)\dnu
=\int\limits_{-\infty}^{\infty}O_r(\nu)\cconjugate{Z_i}(\nu)\dnu\\
 \int\limits_{-\infty}^{\infty}z_r(t)o_i(t)\dt &
=\int\limits_{-\infty}^{\infty}Z_r(\nu)\cconjugate{O_i}(\nu)\dnu
=\int\limits_{-\infty}^{\infty}Z_r(\nu)\cconjugate{O_i}(\nu)\dnu
\end{align*}
Now, using the expressions for $O_r(\nu)$ and $O_i(\nu)$ given above, we have that:
\begin{align*}
              \int\limits_{-\infty}^{\infty}z_i(t)o_r(t)\dt &=
\phantom{-}\xi\int\limits_{-\infty}^{\infty}\sign(\nu)Z_i(\nu)\cconjugate{Z_r}(\nu)\dnu\\
              \int\limits_{-\infty}^{\infty}z_r(t)o_i(t)\dt &=
          -\xi\int\limits_{-\infty}^{\infty}\sign(\nu)Z_r(\nu)\cconjugate{Z_i}(\nu)\dnu
\end{align*}
But, remember that $z_r(t)$, $o_r(t)$, $o_i(t)$ and $z_i(t)$ are real valued, and thus
$\int\limits_{-\infty}^{\infty}z_i(t)o_r(t)\dt\in\R$ and $\int\limits_{-\infty}^{\infty}z_r(t)o_i(t)\dt\in\R$.
As a consequence, the two integrals equal to $\beta$ may be written:
\begin{align*}
\int\limits_{-\infty}^{\infty}z_i(t)o_r(t)\dt &=
          -\Im\left[\int\limits_{-\infty}^{\infty}\sign(\nu)Z_i(\nu)\cconjugate{Z_r}(\nu)\dnu\right]\\
\int\limits_{-\infty}^{\infty}z_r(t)o_i(t)\dt &=
\phantom{-}\Im\left[\int\limits_{-\infty}^{\infty}\sign(\nu)Z_r(\nu)\cconjugate{Z_i}(\nu)\dnu\right]
\end{align*}
As these integrals are all equal, we simply need to show that one of
them is not equal to zero for our purpose. 
Now, let us consider the following quantity:
$$
W(\nu)=\Im\left[\int\limits_{-\infty}^{\infty}Z_r(\nu)\cconjugate{Z_i}(\nu)d\nu\right]
$$
One can see by direct calculation (using Hermitian symmetry of the
classical Fourier transform) that:
\begin{align*}
W(-\nu) &=\phantom{-}\Im\left[\int\limits_{-\infty}^{\infty}Z_r(-\nu)\cconjugate{Z_i}(-\nu)\dnu\right]\\
        &=\phantom{-}\Im\left[\int\limits_{-\infty}^{\infty}\cconjugate{Z_r}(\nu)Z_i(\nu)\dnu\right]\\
        &=         - \Im\left[\int\limits_{-\infty}^{\infty}Z_r(\nu)\cconjugate{Z_i}(\nu)\dnu\right]\\
        &=-W(\nu) 
\end{align*}
which shows that $W(\nu)$ is odd. Now, as a consequence, the following is true:
$$\Im\left[\int\limits_{-\infty}^{\infty}\sign(\nu)Z_r(\nu)\cconjugate{Z_i}(\nu)d\nu\right]
\text{is even} $$ 
And this makes it impossible to cancel out, except in the very
special case where $Z_r(\nu)$ and/or $Z_i(\nu)$ are equal to zero for any $\nu$. 
Thus the following is true:
$$
\Im\left[\int\limits_{-\infty}^{\infty}\sign(\nu)Z_r(\nu)\cconjugate{Z_i}(\nu)d\nu\right]=\int\limits_{-\infty}^{\infty}z_r(t)o_i(t)dt=\beta
\neq 0
$$
This completes the proof.
\end{IEEEproof}

In the most `classical' case ($z(t)$ taking values in $\C^{\i}$), this
means that if one computes the corresponding hyperanalytic signal
using a quaternion Fourier transform with axis \j, then the
hyperanalytic signal $h(t)$ will be $\C^{\j}$-proper.

A way to figure out the meaning of the $\C$-properness of the
hyperanalytic signal  is to see it as a pair of complex signals,
both \emph{improper} and orthogonal to each other
(in the sense of the complex scalar product) but not jointly proper.

We now consider the case when the original signal is already a complex
analytic signal. This is a degenerate case, and the consequences for
properness properties are given as follows.

\begin{property}
\label{prop:properH}
If $z(t)$ with values in $\C^{\mu}$ is a proper complex signal, then
its hyperanalytic signal $h(t)$, computed using a quaternion Fourier
transform with axis $\mu$ is $\C^{\mu}$-superproper\footnote{We use the term
  $\C^{\mu}$-superproper for a $\C^{\mu}$-proper quaternion valued
  signal obtained from a complex proper signal, {\em i.e.} for which
  $\gamma=0$ in the covariance matrix.}, provided that $(1,\xi,\mu,\xi\mu)$ is a
quaternion basis.
 As a consequence, the covariance of $h(t)$ reads:
$$
\Gamma_{h}=
\begin{bmatrix}
\alpha &      0 &      0 & \beta\\
     0 & \alpha &  \beta & 0\\ 
     0 &  \beta & \alpha & 0\\
 \beta &      0 &      0 & \alpha
\end{bmatrix}
$$
with:
$$\alpha=\int\limits_{-\infty}^{\infty}z^2_r(t)\dt=
         \int\limits_{-\infty}^{\infty}o^2_r(t)\dt=
         \int\limits_{-\infty}^{\infty}z^2_i(t)\dt=
         \int\limits_{-\infty}^{\infty}o^2_i(t)\dt
$$
and $\beta$ being as defined in Property \ref{prop:properC}.
\end{property}
\begin{IEEEproof}
Referring to Property \ref{prop:properC} we need to prove that: $\alpha=\omega$ and $\gamma=0$.

First, as $z(t)$ is a complex analytic signal, it is \emph{proper} and so:
$$
\int\limits_{-\infty}^{\infty}z_r(t)z_i(t)\dt=\gamma=0
$$
Then, it is also a well-known property of analytic signals that their
real and imaginary parts have the same magnitude.
Here, this means that:
$$
\int\limits_{-\infty}^{\infty}z^2_r(t)\dt=\int\limits_{-\infty}^{\infty}z^2_i(t)\dt
$$
and consequently $\alpha=\omega$. This completes the proof.

Note that we have only highlighted here the properness properties of
the hyperanalytic signal and left for future work the possiblity of
taking advantage of these properties in its processing.
\end{IEEEproof}

\section{Examples}
\label{Examplessection}

\label{sec:example}
We present here an example signal, with its complex envelope
and phase computed using a hyperanalytic signal.
The hyperanalytic signal is computed from a one-sided quaternion
Fourier transform, with the negative frequencies suppressed
using the discrete-time scheme first published by Marple
\cite{10.1109/78.782222}.

\begin{figure}
\includegraphics[width=\columnwidth]{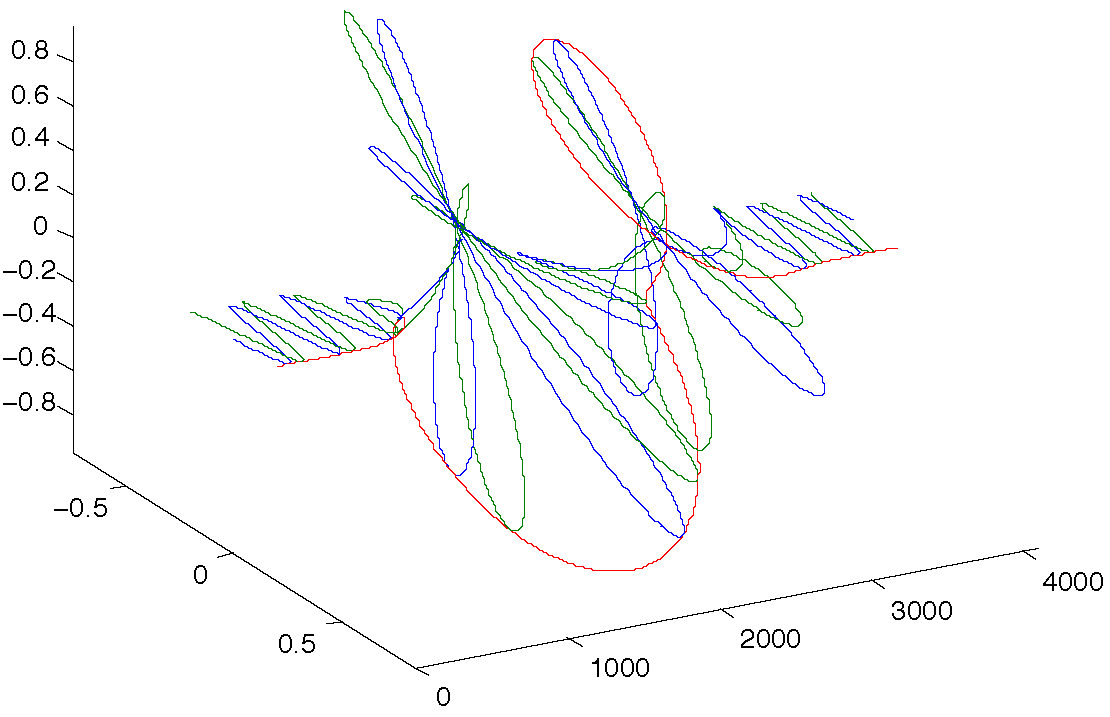}
\caption{\label{fig:helix_gaussian}Sinusoidal carrier modulated by a
complex helix and two Gaussian pulses.
Blue: original signal; green: orthogonal signal; red: complex envelope.}
\end{figure}

\begin{figure}
\includegraphics[width=\columnwidth]{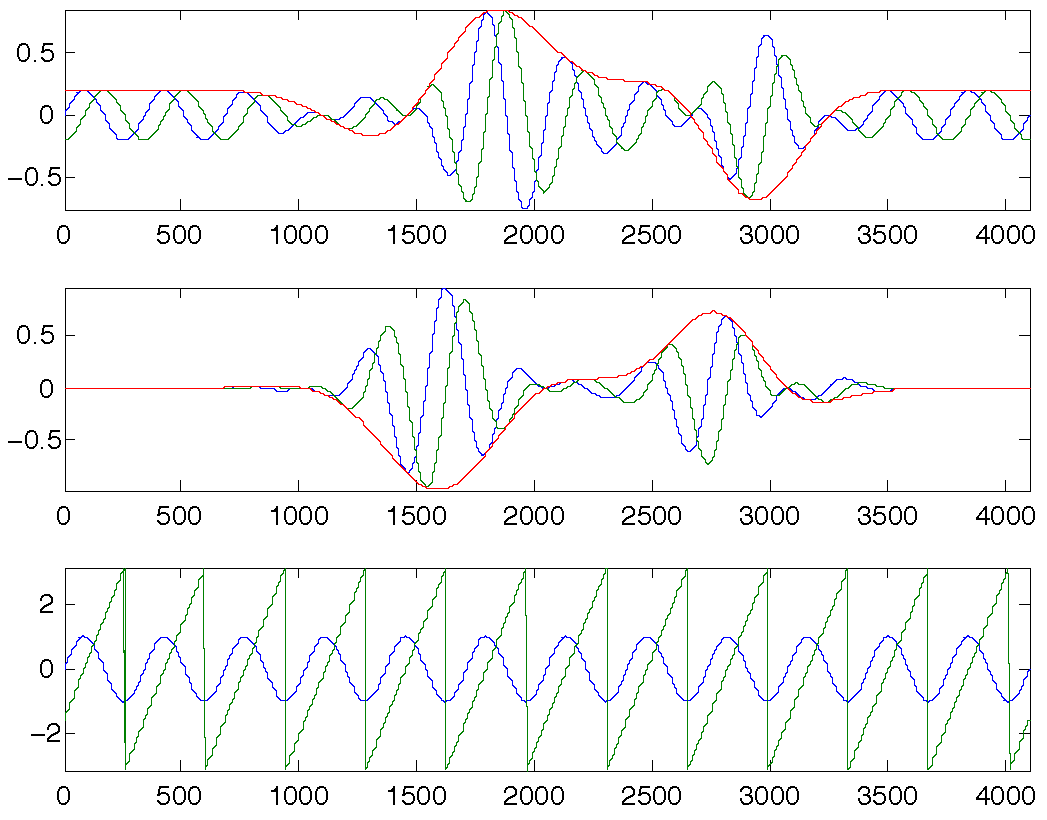}
\caption{\label{fig:helixcomponents}
Top:    real      parts of signals in Figure \ref{fig:helix_gaussian};
Middle: imaginary parts of signals in Figure \ref{fig:helix_gaussian};
Lower:  sinusoidal carrier (blue) and recovered carrier phase $\phi(t)$ (green),
        using Definition \ref{def:cdenv}.}
\end{figure}

Figure \ref{fig:helix_gaussian} shows a complex analogue of
the signal shown in Figure \ref{fig:analytic}.
The sinusoidal carrier is amplitude modulated by two 
Gaussian pulses, and modulated in angle (in the complex plane)
by multiplication with a complex exponential.
The figure shows an orthogonal signal and complex envelope,
computed from a hyperanalytic signal constructed using a
quaternion Fourier transform, with the envelope extracted
as described in Definition \ref{def:cdenv}.
Figure \ref{fig:helixcomponents} shows the real and imaginary components
of the signals shown in Figure \ref{fig:helix_gaussian},
as well as the original sinusoidal carrier and the extracted complex phase
$\phi(t)$ (in this case the imaginary part is negligible, as expected
from Property \ref{prop:instphase} and is not shown).
The envelope can be seen to pass through all four quadrants.
Phase unwrapping was used to eliminate discontinuities in the
complex envelope, operating on the angle of the complex envelope
and constructing the final result from the modulus of the complex
envelope combined with the phase-unwrapped angle.
The QTFM toolbox \cite{qtfm} was used for the computations, and
in particular, the function \texttt{cdpolar.m} was used to extract
the complex envelope from the hyperanalytic signal.

\section{Conclusions}
\label{Conclusionsection}

We have shown in this paper how the classical analytic signal
concept may be extended to the case of the hyperanalytic
signal of an original complex signal, constructing an
orthogonal complex signal by one of three equivalent methods.
The quaternion method yields an interpretation of the
hyperanalytic signal as a quaternion signal which leads
naturally to the definition of the complex envelope.
We have shown that the quaternion polar representation
developed by the authors in 2008 permits us to recover
the phase of a sinusoidal carrier modulated by a complex
signal with frequency content lower than the carrier
frequency, extending the classical result of amplitude
modulation.

We have presented the properness of the hyperanalytic
signal, showing that for a given improper complex signal, its
hyperanalytic signal is $\C$-proper.

As we stated earlier in the paper, our aim in this work
has been to make some progress towards an analytic signal
concept applicable to vector or quaternion signals.
We believe the present paper indicates some of the
problems to be overcome in doing this, including the
need to work in a higher dimensional algebra than the
quaternions, just as the results in this paper are most
simply expressed using quaternion concepts rather than
the concept of complex pairs.

There is also scope to combine the ideas presented here
with the ideas of monogenic signals (that is analytic
signals in two dimensions, or the analytic signal of
images) \cite{FelsbergSommer:2001}.

\bibliographystyle{IEEEtran}
\bibliography{IEEEabrv,hyperanalytic}

\end{document}